\documentclass[12pt,twoside,leqno,openany]{amsart}

\usepackage{amssymb,amsbsy,amsmath,amsfonts,amssymb,amscd,times,
graphics,color,xypic,footmisc,fancyhdr,multicol,fancybox,
graphicx,mathrsfs,rotating,ifthen,wasysym}
\usepackage[all]{xy}

\usepackage[T1]{fontenc}
\newcommand{\C}{\mathbb{C}}

\newcommand{\R}{\mathbb{R}}

\newcommand{\red}{\textcolor{red}}

\newcommand{\HEAD}[2]{%
\pagestyle{fancy}
\fancyhead[RO]{\tiny\sf\thepage}
\fancyhead[CO]{{\tiny\sf #1}}
\fancyhead[LE]{\tiny\sf\thepage}
\fancyhead[CE]{{\tiny\sf #2}}
\fancyfoot{}}

\newtheorem{The}{Theorem}[section]

\newtheorem{Theorem}{Theorem}[section]

\newtheorem{Proposition}[The]{Proposition}
\newtheorem{Lemma}[The]{Lemma}

\theoremstyle{definition}

\renewcommand{\dim}{\text{\footnotesize\sf dim}}

\renewcommand{\exp}{\text{\footnotesize\sf exp}}

\newcommand{\isqrt}{{\scriptstyle{\sqrt{-1}}}}

\renewcommand{\lim}{\text{\footnotesize\sf lim}}

\newcommand{\rank}{\text{\footnotesize\sf rank}}

\newcommand{\smallbullet}{{\scriptscriptstyle{\bullet}}}

\newcommand{\vf}{\vfill\end{document}}

\newcommand{\zero}[1]{\underline{#1}_{\red{\circ}}}

\let\mathcal\mathscr

\begin{document}

$\:$

\bigskip\bigskip

\begin{center}

{\Large\bf Canonical Cartan Connections}

\bigskip

{\Large\bf on Maximally Minimal}

\medskip

{\Large\bf Generic Submanifolds $M^5 \subset \C^4$}

\end{center}

\medskip

\begin{center}
Jo\"el {\sc Merker}, 
Samuel {\sc Pocchiola},
Masoud {\sc Sabzevari}
\end{center}

\bigskip

\begin{center}
\begin{minipage}[t]{10.25cm}
\baselineskip=0.32cm 
{\scriptsize
{\bf Abstract.}
On a real analytic $5$-dimensional CR-generic submanifold 
$M^5 \subset \C^4$ of codimension $3$ hence of CR dimension $1$,
which enjoys the generically satisfied nondegeneracy condition:
\[
\aligned
{\bf 5}
\,=\,
{\sf rank}_\C
\big(
&
T^{1,0}M+T^{0,1}M
+
\big[T^{1,0}M,\,T^{0,1}M\big]
\,+
\\
&
+
\big[T^{1,0}M,\,[T^{1,0}M,T^{0,1}M]\big]
+
\big[T^{0,1}M,\,[T^{1,0}M,T^{0,1}M]\big]
\big),
\endaligned
\]
a canonical Cartan connection is constructed after reduction
to a certain partially explicit $\{ e\}$-structure
of the concerned local biholomorphic equivalence problem.

}
\end{minipage}
\end{center}

\bigskip

\begin{center}
\begin{minipage}[t]{11.75cm}
\baselineskip =0.35cm {\scriptsize

\centerline{\bf Table of contents}

\medskip

{\bf \ref{introduction}.~Introduction
\dotfill~\pageref{introduction}.}

{\bf \ref{initial-Lie-Darboux}.~Initial 
$G$-structure and initial Darboux structure
\dotfill~\pageref{initial-Lie-Darboux}.}

{\bf \ref{reduction-e-structure}.~Reduction of the initial 
$G$-structure to an $\{ e\}$-structure
\dotfill~\pageref{reduction-e-structure}.}

{\bf \ref{Cartan-connection}.~Construction of a Cartan connection
\dotfill~\pageref{Cartan-connection}.}

{\bf \ref{beyond-connections}.~Beyond Cartan connections
\dotfill~\pageref{beyond-connections}.}

}\end{minipage}
\end{center}

\bigskip

\section{\bf Introduction}
\label{introduction}
\HEAD{\ref{introduction}.~Introduction}{
Jo\"el {\sc Merker}, Samuel {\sc Pocchiola}, Masoud {\sc Sabzevari}}

The goal of this announcement is to study local real analytic ($\mathcal{
C}^\omega$) $5$-dimensional CR-generic submanifolds:
\[
M^5 
\subset 
\C^4
\]
of codimension $3$ hence of CR dimension $1$ that are {\sl maximally
minimal} in the sense that:
\[
\aligned
{\bf 5}
\,=\,
\rank_\C
\Big(
&
T^{1,0}M+T^{0,1}M
+
\big[T^{1,0}M,\,T^{0,1}M\big]
\,+
\\
&
+
\big[T^{1,0}M,\,[T^{1,0}M,T^{0,1}M]\big]
+
\big[T^{0,1}M,\,[T^{1,0}M,T^{0,1}M]\big]
\Big).
\endaligned
\] 
In the terminology of~\cite{ Merker-Pocchiola-Sabzevari-5-CR-II},
such CR manifolds $M^5 \subset \C^4$ are said to belong to the 
{\sl General Class
$\text{\sf III}_{\text{\sf 1}}$}. Most considerations being local,
by convention, neighborhoods and their shrinkings will be unmentioned,
as in \'Elie Cartan's original works, {\em cf.} also Peter Olver's
monograph~\cite{ Olver-1995}.

Therefore, if $\mathcal{ L}$ denotes a local vector field generator for 
$T^{1, 0}M$, then the $5$ vector fields:
\[
\mathcal{L},\ \ \ \ 
\overline{\mathcal{L}},\ \ \ \ \
\mathcal{T}:=\isqrt\,\big[\mathcal{L},\overline{\mathcal{L}}\big]
=
\overline{\mathcal{T}},\ \ \ \
\mathcal{S}:=\big[\mathcal{L},\mathcal{T}\big],\ \ \ \ \
\overline{\mathcal{S}}=\big[\overline{\mathcal{L}},\mathcal{T}\big],
\]
are assumed to 
constitute a (local) frame for $\C \otimes_\R TM$, that is to say, at
each point $p \in M$:
\[
{\bf 5}
=
\rank_\C
\Big(
\mathcal{L}\big\vert_p,\,\,
\overline{\mathcal{L}}\big\vert_p,\,\,
\mathcal{T}\big\vert_p,\,\,
\mathcal{S}\big\vert_p,\,\,
\overline{\mathcal{S}}\big\vert_p
\Big).
\]

In coordinates $(z, w_1, w_2, w_3) \in \C^4$ with $w_j = u_j + 
\isqrt\, v_j$, for Beloshapka's {\sf c}ubic model having
equations:
\[
M_{\sf c}^5\colon
\ \ \ \ \
\left[
\aligned
v_1
&
=
z\overline{z},
\\
v_2
&
=
z^2\overline{z}+z\overline{z}^2,
\\
v_3
&
=
-\,\isqrt\,\big(
z^2\overline{z}
-
z\overline{z}^2
\big),
\endaligned\right.
\]
the five vector fields in question (\cite{ Merker-Sabzevari-III-1},
p.~94) are visibly everywhere linearly independent:
\[
\aligned
\mathcal{L}_{\sf c}
&
=
{\textstyle{\frac{\partial}{\partial z}}}
+
\isqrt\,\overline{z}\,
{\textstyle{\frac{\partial}{\partial u_1}}}
+
\isqrt\,\big(
2\,z\overline{z}+\overline{z}^2\big)\,
{\textstyle{\frac{\partial}{\partial u_2}}}
+
\big(2\,z\overline{z}-\overline{z}^2\big)\,
{\textstyle{\frac{\partial}{\partial u_3}}},
\\
\overline{\mathcal{L}}_{\sf c}
&
=
{\textstyle{\frac{\partial}{\partial\overline{z}}}}
-
\isqrt\,z\,
{\textstyle{\frac{\partial}{\partial u_1}}}
-
\isqrt\,
\big(
2z\overline{z}+z^2
\big)\,
{\textstyle{\frac{\partial}{\partial u_2}}}
+
\big(
2\,z\overline{z}
-
z^2
\big)\,
{\textstyle{\frac{\partial}{\partial u_3}}},
\\
\mathcal{T}_{\sf c}
&
=
\ \ \ \ \ \ \ \ \ \ \ \ \ \ \ \
2\,
{\textstyle{\frac{\partial}{\partial u_1}}}
+
4\,\big(z+\overline{z}\big)\,
{\textstyle{\frac{\partial}{\partial u_2}}}
-
4\,\isqrt\,
\big(z-\overline{z}\big)\,
{\textstyle{\frac{\partial}{\partial u_3}}},
\\
\mathcal{S}_{\sf c}
&
=
\ \ \ \ \ \ \ \ \ \ \ \ \ \ \ \ \ \ \ \ \ \ \ \ \ \ \ \ \ \ \ 
\ \ \ \ \ \ \ \ \ \ \ \
4\,
{\textstyle{\frac{\partial}{\partial u_2}}}
-
4\,\isqrt\,
{\textstyle{\frac{\partial}{\partial u_3}}},
\\
\overline{\mathcal{S}}_{\sf c}
&
=
\ \ \ \ \ \ \ \ \ \ \ \ \ \ \ \ \ \ \ \ \ \ \ \ \ \ \ \ \ \ \ 
\ \ \ \ \ \ \ \ \ \ \ \
4\,
{\textstyle{\frac{\partial}{\partial u_2}}}
+
4\,\isqrt\,
{\textstyle{\frac{\partial}{\partial u_3}}}.
\endaligned
\]

Inspired by other similar models such as the Heisenberg {\sf q}uadric
$M_{\sf q}^3 \subset \C^2 \ni (z, w_1)$ the equation of which is the 
just the first one $v_1 = z\overline{ z}$ above, inspired also
by the companion Beloshapka cubic $M_{\sf c}^4 \subset \C^3 \ni
(z, w_1, w_2)$ the two equations of which are just the first
two above $v_1 = z\overline{ z}$, $v_2 = z^2\overline{ z} + 
z\overline{ z}^2$, knowing that canonical Cartan connections have been
constructed for general geometry-preserving deformations of these two
models, by, respectively, Cartan~\cite{ Cartan-1932} and by
Beloshapka-Ezhov-Schmalz~\cite{ BES-2007}, the objective of the
present announcement is to show how to extract, from the recent
prepublication of an extensive memoir (\cite{ Merker-Sabzevari-III-1}),
the construction of a canonical Cartan connection associated to
local biholomorphic equivalences\,\,---\,\,or equivalently, to real
analytic CR equivalences\,\,---\,\,for General Class $\text{\sf
III}_{\text{\sf 1}}$ CR manifolds.

Beloshapka's cubic model $M_{\sf c}^5 \subset \C^4$
happens to be a {\sl homogeneous space}, namely a quotient:
\[
M_{\sf c}^5
\cong
G^7
\big/
H^2
\cong
N_4^5,
\]
of a $7$-dimensional real Lie group $G^7$ by a $2$-dimensional closed
commutative Lie subgroup $N^2 \cong (\C^*, \times)$, the resulting quotient
being the unique connected and simply connected {\em nilpotent}
Lie group corresponding to the real {\em nilpotent} Lie
algebra with generators ${\sf x}_1$, ${\sf x}_2$, ${\sf x}_3$,
${\sf x}_4$, ${\sf x}_5$ named $\mathfrak{ n}_5^4$
in the Goze-Remm classification:
\[
\mathfrak{n}_5^4
\colon
\ \ \ \ \ \ \ \ \ \
\left\{
\aligned
{}
[{\sf x}_1,{\sf x}_2]
&
=
{\sf x}_3,
\\
[{\sf x}_1,{\sf x}_3]
&
=
{\sf x}_4,
\\
[{\sf x}_2,{\sf x}_3]
&
=
{\sf x}_5,
\endaligned\right.
\]
unwritten brackets being zero.

Indeed, it is known ({\em see}~\cite{ Merker-Sabzevari-III-1} for
details) that the Lie algebra 
$\mathfrak{ aut}_{\text{\sc cr}} ( M_{\sf c}^5)$ 
of infinitesimal CR automorphisms of 
Beloshapka's
cubic model $M_{\sf c}^5 \subset \C^4$ is $7$-dimensional, generated
by the real parts of the following $7$ vector fields of type $(1, 0)$
with holomorphic coefficients:
\[
\aligned 
S_2
& 
:=
{\textstyle{\frac{\partial}{\partial w_3}}},
\\
S_1 
& 
:=
{\textstyle{\frac{\partial}{\partial w_2}}},
\\
T
& 
:=
{\textstyle{\frac{\partial}{\partial w_1}}},
\\
L_1 
& 
:=
{\textstyle{\frac{\partial}{\partial z}}}
+
2\,\isqrt\,z\,
{\textstyle{\frac{\partial}{\partial w_1}}}
+
\big(2\,\isqrt\,z^2+4\,w_1\big)\,
{\textstyle{\frac{\partial}{\partial w_2}}}
+
2\,z^2\,
{\textstyle{\frac{\partial}{\partial w_3}}},
\\
L_2
&
:=
\isqrt\,
{\textstyle{\frac{\partial}{\partial z}}}
+
2\,z\,
{\textstyle{\frac{\partial}{\partial w_1}}}
+
2\,z^2\,
{\textstyle{\frac{\partial}{\partial w_2}}}
-
\big(2\,\isqrt\,z^2-4\,w_1\big)\,
{\textstyle{\frac{\partial}{\partial w_3}}},
\\
D
&
:=
z\,
{\textstyle{\frac{\partial}{\partial z}}}
+
2\,w_1\,
{\textstyle{\frac{\partial}{\partial w_1}}}
+
3\,w_2\,
{\textstyle{\frac{\partial}{\partial w_2}}}
+
3\,w_3\,
{\textstyle{\frac{\partial}{\partial w_3}}},
\\
R
&
:=
\isqrt\,z\,
{\textstyle{\frac{\partial}{\partial z}}}
-
w_3\,
{\textstyle{\frac{\partial}{\partial w_2}}}
+
w_2\,
{\textstyle{\frac{\partial}{\partial w_3}}}.
\endaligned
\]
Their expressions show well that 
the isotropy Lie subalgebra of the origin $0 \in \C^4$ is
generated by the last two fields $D$ and $R$, the only ones all
of whose coefficients vanish there. Furthermore, the complete 
Lie bracket commutator table:

\begin{center}
\begin{tabular} [t] { c | c c c c c c c }
& $S_2$ & $S_1$ & $T$ & $L_2$ & $L_1$ & $D$ & $R$
\\
\hline $S_2$ & $0$ & $0$ & $0$ & $0$ & $0$ & $3S_2$ & $-S_1$
\\
$S_1$ & $*$ & $0$ & $0$ & $0$ & $0$ & $3S_1$ & $S_2$
\\
$T$ & $*$ & $*$ & $0$ & $4S_2$ & $4S_1$ & $2T$ & $0$
\\
$L_2$ & $*$ & $*$ & $*$ & $0$ & $-4T$ & $L_2$ & $-L_1$
\\
$L_1$ & $*$ & $*$ & $*$ & $*$ & $0$ & $L_1$ & $L_2$
\\
$D$ & $*$ & $*$ & $*$ & $*$ & $*$ & $0$ & $0$
\\
$R$ & $*$ & $*$ & $*$ & $*$ & $*$ & $*$ & $0$
\end{tabular}
\end{center}

\medskip\noindent
shows that the five fields $S_2$, $S_1$, $T$, $L_2$, $L_1$ generate a
nilpotent Lie subalgebra of $\mathfrak{ aut}_{\text{\sc cr}} 
( M_{\sf c}^5)$ visibly isomorphic to $\mathfrak{ n}_5^4$.

Once an expected appropriate model geometry $G / H$ of Klein type 
has been set up, its curved Cartan type deformations can enter the scene.

Recall that a {\sl Cartan geometry} on a $\mathcal{ C}^\omega$
manifold $M$ {\sl modelled} on a homogeneous space $G/H$, where $G$ is
a connected Lie group and $H \subset G$ is a closed connected Lie
subgroup, with ${\rm Lie}(H) = \mathfrak{h} \,\subset\, \mathfrak{g} =
{\rm Lie}(G)$, is a pair $(P,\, \varpi)$ consisting of an
$H$-principal bundle $\pi \colon P \longrightarrow M$, with right
action $R_h \colon p \longmapsto p\, h$ for $h \in H$ and $p \in P$,
together with a $\mathfrak{ g}$-valued differential $1$-form $\varpi
\colon TP \,\longrightarrow\, \mathfrak{g}$ enjoying:

\begin{itemize}

\smallskip\item[{\bf (i)}]\,
$\varpi_p \colon T_p P \overset{\sim}{\longrightarrow} \mathfrak{ g}$ 
is isomorphic a every point $p \in P$; 

\smallskip\item[{\bf (ii)}]\,
for every ${\sf y} \in \mathfrak{ h}$, if $Y^\dag$ denotes the vector
field $Y^\dag \vert_p := 
\frac{ d}{ dt} \big\vert_0 \big( p \, \exp( t\, {\sf y}) \big)$, then
$\varpi ( Y^\dag) = {\sf y}$; 

\smallskip\item[{\bf (iii)}]\,
at every $p \in P$, for every $v_p \in T_p P$:
\[
\varpi_{ph}
\big(
R_{h*}(v_p)
\big)
\,=\,
{\rm Ad}(h^{-1})\,
\big[
\varpi_p(v_p)
\big].
\] 

\end{itemize}

If $M$ is a manifold endowed with a certain determined type of
geometric structure, {\em e.g.} an integrable CR structure, either
abstract or embedded, it is adequate to call {\sl canonical} a Cartan
geometry $(P, \varpi)$ on $M$ when all (local or global) automorphisms
${\sf H} \colon\, M \longrightarrow M$ of the geometric structure {\em
lift} as bundle automorphisms:
\[
\xymatrix{
P \ar[r]^{\widehat{{\sf H}}} \ar[d]_\pi
& 
P \ar[d]^\pi
\\
M \ar[r]_{\sf H} & M,
}
\]
satisfying $\widehat{\sf H}^*\big( \varpi\big) = \varpi$, and when
conversely also, all such $\pi$-fiber preserving maps $\widehat{\sf
H}$ satisfying $\widehat{\sf H}^*\big( \varpi \big) = \varpi$, {\em
descend} as automorphisms ${\sf H} \colon\, M \to M$.

The main result here, computationally less advanced than what comes out 
from~\cite{ Merker-Sabzevari-III-1}, can be summarized as follows.

\begin{Theorem}
\label{main-theorem}
Associated to every Class $\text{\sf III}_{\text{\sf 1}}$ local
CR-generic $\mathcal{ C}^\omega$ submanifold $M^5 \subset \C^4$, there
is a canonical Cartan connection $(P^7, \varpi)$ modelled on the
nilpotent homogeneous space $G^7 \big/ H^2 \cong N_4^5
\cong M_{\sf c}^5$ whose natural orbit
space is Beloshapka's cubic model $M_{\sf c}^5 \subset \C^4$.
\end{Theorem}

Generally, if $(P, \varpi)$ is a canonical Cartan
connection on a manifold $M$ belonging to a determined general class
of geometric structures, then it automatically carries an associated
canonical {\sl absolute parallelism}, or so-called {\sl $\{
e\}$-structure}, obtained, with a basis ${\sf e}_1, \dots, {\sf
e}_r$ of $\mathfrak{ g}$ and with $r := \dim_\R \mathfrak{ g}$, by
plainly setting:
\[
V_i\big\vert_p
\,:=\,
\varpi_p^{-1}({\sf e}_i)
\ \ \ \ \ \ \ \ \ \ \ \ \
{\scriptstyle{(1\,\leqslant\,i\,\leqslant\,r)}},
\]
the obtained vector fields $V_1, \dots, V_r$ making up a frame on
$P$. 
 
The proof of the theorem consists in constructing an absolute
parallelism (Section~\ref{reduction-e-structure}) 
and to verify that it satisfies
the algebraic conditions required to constitute a Cartan connection
(Section~\ref{Cartan-connection}).

\medskip\noindent{\bf Acknowledgments.} The authors would like to
thank Professor Alexander Isaev and an anonymous referee, for
realizing that the plain constructions of absolute parallelisms or of
Cartan connections usually stays at a lower level than precise
inspections of the parametric expressions of incoming invariants and
of their {\em nonlinear} relations.

\section{\bf Initial $G$-structure and initial Darboux structure}
\label{initial-Lie-Darboux}
\HEAD{\ref{initial-Lie-Darboux}.~Initial 
$G$-structure and initial Darboux structure}{
Jo\"el {\sc Merker}, Samuel {\sc Pocchiola}, Masoud {\sc Sabzevari}}

The goal being to set up a Cartan procedure in order to reduce to an
$\{ e\}$-structure the equivalence problem under local biholomorphisms
for such $M^5 \subset \C^4$ belonging to Class $\text{\sf
III}_{\text{\sf 1}}$, and ultimately, to construct an associated
Cartan connection, the first task is to examine how such a frame
$\big\{ \mathcal{ L}, \overline{ \mathcal{ L}}, \mathcal{ T},
\mathcal{ S}, \overline{ \mathcal{ S}} \big\}$ transforms under an
arbitrary local biholomorphic map:
\[
{\sf H}\colon\ \ \
\C^4
\,\longrightarrow\,
{\C'}^4.
\]
After appropriate unmentioned shrinkings of concerned neighborhoods,
the image:
\[
{\sf H}(M)
=:
M'
\]
also becomes (\cite{ Merker-Pocchiola-Sabzevari-5-CR-II}, Section~3) 
a CR-generic $3$-codimensional submanifold ${M'}^5 \subset \C^4$
of CR dimension $1$, also equipped from its side with an analogous,
independently constructed, frame starting from
a local generator $\mathcal{ L}'$ of $T^{1, 0} M'$:
\[
\big\{
\mathcal{L}',\,
\overline{\mathcal{L}}',\,
\mathcal{T}',\,
\mathcal{S}',\,
\overline{\mathcal{S}}'
\big\},
\ \ \ \ \
\text{\rm with}\ \
\mathcal{T}'
:=
\isqrt\,\big[\mathcal{L}',\overline{\mathcal{L}}'\big],\ \
\mathcal{S}'
:=
\big[\mathcal{L}',\mathcal{T}'\big],
\]
and because the restriction 
\[
{\sf H}
\big\vert_M\colon M 
\,\longrightarrow\,
M'
\]
is known to be a CR-diffeomorphism, namely because
${\sf H}_* ( T^{1, 0} M ) = T^{1, 0} M'$, there must exist a
nowhere vanishing $\mathcal{ C}^\omega$ function defined
on $M'$ such that (\cite{ Merker-Pocchiola-Sabzevari-5-CR-II}, 
Section~4):
\[
{\sf H}_*(\mathcal{L})
=
a'\,\mathcal{L}',
\]
whence by plain conjugation (conventionally not bearing on the
differential ${\sf H}_*$):
\[
{\sf H}_*
\big(\overline{\mathcal{L}}\big)
=
\overline{a}'\,
\overline{\mathcal{L}}'.
\]

\begin{Proposition}
\label{transfer-frame}
{\rm (\cite{ Merker-Sabzevari-III-1})}
Under a local biholomorphic map ${\sf H} \colon \C^4 \to {\C'}^4$, 
two adapted frames associated to two ${\sf H}$-equivalent CR-generic 
$M^5 \subset \C^4$ and ${\sf H}(M) =: {M'}^5 \subset {\C'}^4$ transfer, 
in terms of certain five $\C$-valued local $\mathcal{ C}^\omega$ 
functions $a'$, $b'$, $c'$, $d'$, $e'$, defined on $M'$ as:
\[
\left(\!\!
\begin{array}{c}
\mathcal{L}
\\
\overline{\mathcal{L}}
\\
\mathcal{T}
\\
\mathcal{S}
\\
\overline{\mathcal{S}}
\end{array}
\!\!\right)
\,=\,
\left(\!\!
\begin{array}{ccccc}
a' & 0 & 0 & 0 & 0
\\
0 & \overline{a}' & 0 & 0 & 0
\\
b' & \overline{b}' & a'\overline{a}' & 0 & 0
\\
e' & d' & c' & a'a'\overline{a}' & 0
\\
\overline{d}' & \overline{e}' & \overline{c}' & 0 &
a'\overline{a}'\overline{a}'
\end{array}
\!\!\right)
\left(\!\!
\begin{array}{c}
\mathcal{L}'
\\
\overline{\mathcal{L}}'
\\
\mathcal{T}'
\\
\mathcal{S}'
\\
\overline{\mathcal{S}}'
\end{array}
\!\!\right).
\]
\end{Proposition}

\proof
Computing further the bracket shows:
\[
\aligned
{\sf H}_*\big(\mathcal{T}\big)
&
=
{\sf H}_*\big(\isqrt\,\big[\mathcal{L},\,\overline{\mathcal{L}}\big]\big)
\\
&
=
\isqrt\,
\big[
{\sf H}_*(\mathcal{L}),\,
{\sf H}_*(\overline{\mathcal{L}})
\big]
\\
&
=
\isqrt\,
\big[
a'\,\mathcal{L}',\,
\overline{a}'\,\overline{\mathcal{L}}'
\big]
\\
&
=
a'\overline{a}'\,
\isqrt\,
\big[\mathcal{L}',\,\overline{\mathcal{L}}'\big]
\underbrace{
-\isqrt\,\overline{a}'\,\overline{\mathcal{L}}'(a')}_{
=:\,b'}
\cdot
\mathcal{L}'
+
\isqrt\,a'\,\mathcal{L}'(\overline{a}')
\cdot
\overline{\mathcal{L}}'
\\
&
=:
a'\overline{a}'\,\mathcal{T}'
+
b'\,\mathcal{L}'
+
\overline{b}'\,\overline{\mathcal{L}}',
\endaligned
\]
in terms of some new function $b'$ to which an independent
name is given. Then quite similarly
(\cite{ Merker-Sabzevari-III-1}, p.~98):
\[
\aligned
{\sf H}_*(\mathcal{S})
&
=
a'a'\overline{a}'\,\mathcal{S}'
+
c'\,\mathcal{T}'
+
e'\,\mathcal{L}'
+
d'\,\overline{\mathcal{L}}',
\\
{\sf H}_*(\overline{\mathcal{S}})
&
=
a'\overline{a}'\overline{a}'\,
\overline{\mathcal{S}}'
+
\overline{c}'\,\mathcal{T}'
+
\overline{d}'\,\overline{\mathcal{L}}'
+
\overline{e}'\,\mathcal{L}',
\endaligned
\]
which completes the proof.
\endproof

The $10$ Lie brackets between the $5$ local vector fields
$\mathcal{ L}$, $\overline{ \mathcal{ L}}$, $\mathcal{ T}$,
$\mathcal{ S}$, $\overline{\mathcal{ S}}$ involve a first set
of $5$ local real analytic local functions $A$, $B$, $P$, $Q$, $R$
appearing in:
\[
\aligned
\big[\mathcal{L},\mathcal{S}\big]
&
=
P\cdot\mathcal{T}
+
Q\cdot\mathcal{S}
+
R\cdot\overline{\mathcal{S}},
\\
\big[\overline{\mathcal{L}},\mathcal{S}\big]
&
=
A\cdot\mathcal{T}+B\cdot\mathcal{S}+
\overline{B}\cdot\overline{\mathcal{S}},
\\
\big[\mathcal{L},\overline{\mathcal{S}}\big]
&
=
A\cdot\mathcal{T}+B\cdot\mathcal{S}+
\overline{B}\cdot\overline{\mathcal{S}},
\\
\big[\overline{\mathcal{L}},\overline{\mathcal{S}}\big]
&
=
\overline{P}\cdot\mathcal{T}+\overline{R}\cdot\mathcal{S}+
\overline{Q}\cdot\overline{\mathcal{S}},
\endaligned
\]
with $A$ being real-valued, as follows from an inspection of
Jacobi identities (\cite{ Merker-Sabzevari-III-1}, Lemma~13.3).

The $5$ further real analytic functions $E$, $F$, $G$, $J$, $K$,
with $J$ also real-valued, appearing in the remaining $3$ Lie brackets:
\[
\aligned
\big[\mathcal{T},\mathcal{S}\big]
&
=
E\cdot\mathcal{T}
+
F\cdot\mathcal{S}
+
G\cdot\overline{\mathcal{S}},
\\
\big[\mathcal{T},\overline{\mathcal{S}}\big]
&
=
\overline{E}\cdot\mathcal{T}
+
\overline{G}\cdot\mathcal{S}
+
\overline{F}\cdot\overline{\mathcal{S}},
\\
\big[\mathcal{S},\overline{\mathcal{S}}\big]
&
=
\isqrt\,J\cdot\mathcal{T}
+
K\cdot\mathcal{S}
-
\overline{K}\cdot\overline{\mathcal{S}},
\endaligned
\]
all express in terms of $A$, $B$, $P$, $Q$, $R$, for instance:
\[
\aligned
E
&
=
\isqrt\,
\Big(
\mathcal{L}(A)
-
\overline{\mathcal{L}}(P)
-
A\,Q
-
\overline{P}\,R
+
B\,P
+
A\,\overline{B}
\Big),
\\
F
&
=
\isqrt\,
\Big(
\mathcal{L}(B)
-
\overline{\mathcal{L}}(Q)
-
R\,\overline{R}
+
B\,\overline{B}
+
A
\Big),
\\
G
&
=
\isqrt\,
\Big(
\mathcal{L}(\overline{B})
-
\overline{\mathcal{L}}(R)
+
B\,R
+
\overline{B}\,\overline{B}
-
R\,\overline{Q}
-
\overline{B}\,Q
-
P
\Big),
\endaligned
\]
with similar longer expressions for $J$ and $K$ unwritten here
(\cite{ Merker-Sabzevari-III-1}, Lemma~13.5).

\medskip

Introduce then the coframe of $\C$-valued $1$-forms on $M$:
\[
\big\{
\overline{\sigma}_0,\,
\sigma_0,\,
\rho_0,\,
\overline{\zeta}_0,\,
\zeta_0
\big\},
\]
which is {\em dual} to the frame $\big\{ \overline{ \mathcal{ S}},\,
\mathcal{ S},\, \mathcal{T},\, \overline{ \mathcal{L}},\,
\mathcal{ L} \big\}$, in this order. The $10$ two-forms making up
a basis of $\C\otimes_\R \Lambda^2 T^*M$ will be constantly ordered as:
\[
\aligned
\overline{\sigma}_0\wedge\sigma_0,
\ \ \ \ \ \ \ \ \ \
\overline{\sigma}_0\wedge\rho_0,
\ \ \ \ \ \ \ \ \ \
\overline{\sigma}_0\wedge\overline{\zeta}_0,
\ \ \ \ \ \ \ \ \ \
\overline{\sigma}_0\wedge\zeta_0,
\\
\sigma_0\wedge\rho_0,
\ \ \ \ \ \ \ \ \ \
\sigma_0\wedge\overline{\zeta}_0,
\ \ \ \ \ \ \ \ \ \
\sigma_0\wedge\zeta_0,
\\
\rho_0\wedge\overline{\zeta}_0,
\ \ \ \ \ \ \ \ \ \
\rho_0\wedge\zeta_0,
\\
\overline{\zeta}_0\wedge\zeta_0.
\endaligned
\]
Then the above initial
Lie bracket structure translates as the following initial
Darboux structure (\cite{ Merker-Sabzevari-III-1}, p.~92)
for the exterior differentials of the basis initial $1$-forms:
\[
\aligned
d\overline{\sigma}_0
&
=
-\,\overline{K}\cdot
\overline{\sigma}_0\wedge\sigma_0
+
\overline{F}\cdot
\overline{\sigma}_0\wedge\rho_0
+
\overline{Q}\cdot
\overline{\sigma}_0\wedge\overline{\zeta}_0
+
\overline{B}\cdot
\overline{\sigma}_0\wedge\zeta_0
\,+
\\
&
\ \ \ \ \ \ \ \ \ \ \ \ \ \ \ \ \ \ \ \ \ \ \ \ \ \ \ \
+
G\cdot
\sigma_0\wedge\rho_0
+
\overline{B}\cdot
\sigma_0\wedge\overline{\zeta}_0
+
R\cdot
\sigma_0\wedge\zeta_0
\,+
\\
&
\ \ \ \ \ \ \ \ \ \ \ \ \ \ \ \ \ \ \ \ \ \ \ \ \ \ \ \ \ \ \ \ \ \
\ \ \ \ \ \ \ \ \ \ \ \ \ \ \ \ \
+
\rho_0\wedge\overline{\zeta}_0,
\\
d\sigma_0
&
=
K\cdot
\overline{\sigma}_0\wedge\sigma_0
+
\overline{G}\cdot
\overline{\sigma}_0\wedge\rho_0
+
\overline{R}\cdot
\overline{\sigma}_0\wedge\overline{\zeta}_0
+
B\cdot
\overline{\sigma}_0\wedge\zeta_0
\,+
\\
&
\ \ \ \ \ \ \ \ \ \ \ \ \ \ \ \ \ \ \ \ \ \ \ \ 
+
F\cdot
\sigma_0\wedge\rho_0
+
B\cdot\sigma_0\wedge\overline{\zeta}_0
+
Q\cdot
\sigma_0\wedge\zeta_0
\,+
\\
&
\ \ \ \ \ \ \ \ \ \ \ \ \ \ \ \ \ \ \ \ \ \ \ \ \ \ \ \ \ \ \ \ \ \ \
\ \ \ \ \ \ \ \ \ \ \ \ \ \ \ \ \ \ \ \ \ \ \ \ \ \ \ \ \ \ \ \ \ \ \
+
\rho_0\wedge\zeta_0,
\\
d\rho_0
&
=
\isqrt\,J\cdot
\overline{\sigma}_0\wedge\sigma_0
+
\overline{E}\cdot
\overline{\sigma}_0\wedge\rho_0
+
\overline{P}\cdot
\overline{\sigma}_0\wedge\overline{\zeta}_0
+
A\cdot
\overline{\sigma}_0\wedge\zeta_0
\,+
\\
&
\ \ \ \ \ \ \ \ \ \ \ \ \ \ \ \ \ \ \ \ \ \ \ \ \ \ \ \ \ \ \
+
E\cdot
\sigma_0\wedge\rho_0
+
A\cdot
\sigma_0\wedge\overline{\zeta}_0
+
P\cdot
\sigma_0\wedge\zeta_0
-\,
\\
&
\ \ \ \ \ \ \ \ \ \ \ \ \ \ \ \ \ \ \ \ \ \ \ \ \ \ \ \ \ \ \ \ \ \ \
\ \ \ \ \ \ \ \ \ \ \ \ \ \ \ \ \ \ \ \ \ \ \ \ \ \ \ \ \ \ \ \ \ \ \
\ \ \ \ \ \ \ \ 
-
\isqrt\,
\overline{\zeta}_0\wedge\zeta_0,
\endaligned
\]
while $d \zeta_0 = 0$ and $d\overline{ \zeta}_0 = 0$.

\section{\bf Reduction of the initial $G$-structure 
to an $\{ e\}$-structure}
\label{reduction-e-structure}
\HEAD{\ref{reduction-e-structure}.~Reduction of the initial 
$G$-structure to an $\{ e\}$-structure}{
Jo\"el {\sc Merker}, Samuel {\sc Pocchiola}, Masoud {\sc Sabzevari}}

In accordance with Proposition~\ref{transfer-frame}, the initial
$G$-structure encoding local biholomorphic equivalences of manifolds
$M^5 \subset \C^4$ belonging to Class $\text{\sf III}_1$ is, after
reordering the frame as $\big\{ \overline{ \mathcal{ S}}, \mathcal{
S}, \mathcal{ T}, \overline{\mathcal{ L}}, \mathcal{ L} \big\}$, the
following closed Lie subgroup of ${\sf GL}_5 ( \C)$:
\[
{\sf G}_{\text{\sf III}_{\text{\sf 1}}}
\,:=\,
\left\{
g\,:=\,
\left(\!
\begin{array}{ccccc}
{\sf a}\overline{\sf a}\overline{\sf a} & 0 & 
\overline{\sf c} & \overline{\sf c} & \overline{\sf d}
\\
0 & {\sf a}{\sf a}\overline{\sf a} & {\sf c} & {\sf d} & {\sf e}
\\
0 & 0 & {\sf a}\overline{\sf a} & \overline{\sf b} & {\sf b}
\\
0 & 0 & 0 & \overline{\sf a} & 0
\\
0 & 0 & 0 & 0 & {\sf a}
\end{array}
\!\right)
\colon\,\,
{\sf a}\,\in\,\C\backslash\{0\},\,\,
{\sf b},\,{\sf c},\,{\sf d},\,{\sf e}
\,\in\,
\C
\right\}.
\]
Following Cartan and Olver (\cite{ Olver-1995}),
transposing then this matrix in order to express
how {\em co}frames do transfer, introduce the
so-called {\sl lifted coframe}:
\[
\aligned
\left(\!\!
\begin{array}{c}
\overline{\sigma}
\\
\sigma
\\
\rho
\\
\overline{\zeta}
\\
\zeta
\end{array}
\!\!\right)
=
\left(\!\!
\begin{array}{ccccc}
{\sf a}{\sf a}\overline{\sf a} & 0 & 0 & 0 & 0
\\
0 & {\sf a}\overline{\sf a}\overline{\sf a} & 0 & 0 & 0
\\
\overline{\sf c} & {\sf c} & {\sf a}\overline{\sf a} & 0 & 0
\\
\overline{\sf e} & {\sf d} & \overline{\sf b} & \overline{\sf a} & 0
\\
\overline{\sf d} & {\sf e} & {\sf b} & 0 & {\sf a}
\end{array}
\!\!\right)
\left(\!\!
\begin{array}{c}
\overline{\sigma}_0
\\
\sigma_0
\\
\rho_0
\\
\overline{\zeta}_0
\\
\zeta_0
\end{array}
\!\!\right),
\endaligned
\]
with ${\sf a}$, ${\sf b}$, ${\sf c}$, ${\sf d}$, ${\sf e}$ here being
{\em independent} variables, replacing the (unknown) functions from
Proposition~\ref{transfer-frame}. The inverse matrix is:
\[
g^{-1}
=
\left(
\begin{array}{ccccc}
\frac{1}{{\sf a}\overline{\sf a}^2} & 0 & 0 & 0 & 0 
\\
0 & \frac{1}{{\sf a}^2\overline{\sf a}} & 0 & 0 & 0 
\\
-\frac{\overline{\sf c}}{{\sf a}^2\overline{\sf a}^3} & 
-\frac{\sf c}{{\sf a}^3\overline{\sf a}^2} &
\frac{1}{{\sf a}\overline{\sf a}} & 0 & 0 
\\
\frac{\overline{\sf b}\,\overline{\sf c}
- \overline{\sf e}{\sf a}\overline{\sf
a}}{{\sf a}^2\overline{\sf a}^4} 
&
\frac{\overline{\sf b}\sf c-
{\sf a}\overline{a}\sf d}{{\sf a}^3\overline{\sf a}^3}
& -\frac{\overline{\sf b}}{{\sf a}\overline{\sf a}^2} &
\frac{1}{\overline{\sf a}} & 0 
\\
\frac{{\sf b}\overline{\sf c}-
{\sf a}\overline{\sf a}\overline{\sf d}}{{\sf
a}^3\overline{\sf a}^3} &
\frac{{\sf bc}-{\sf e}{\sf a}
\overline{\sf a}}{{\sf a}^4\overline{\sf a}^2} & 
-\frac{\sf b}{{\sf a}^2\overline{\sf a}} & 0 & \frac{1}{{\sf a}} 
\end{array}
\right).
\]

Then the
Maurer-Cartan forms of this matrix group ${\sf G}_{ \text{\sf III}_{
\text{\sf 1}}}$ appear in the full expression of $dg \cdot g^{
-1}$, a new $5 \times 5$ matrix which happens to be of the form:
\[
dg\cdot g^{-1}
\,=\,
\left(\!
\begin{array}{ccccc}
\alpha^1+2\overline{\alpha}^1 & 0 & 0 & 0 & 0
\\
0 & 2\alpha^1+\overline{\alpha}^1 & 0 & 0 & 0
\\
\overline{\alpha}^2 & \alpha^2 & \alpha^1+\overline{\alpha}^1 & 0 & 0
\\
\overline{\alpha}^3 & \overline{\alpha}^4 & \overline{\alpha}^5 &
\overline{\alpha}^1 & 0
\\
\alpha^4 & \alpha^3 & \alpha^5 & 0 & \alpha^1
\end{array}
\!\right),
\]
where (the expressions of $\alpha^3$, $\alpha^4$, $\alpha^5$
will not be used):
\[
\alpha^1
:=
\frac{d{\sf a}}{{\sf a}},
\ \ \ \ \ \ \ \ \ \ \ \ \ \ \ \ \ \ 
\alpha^2
:=
-\,
\frac{{\sf c}\,d{\sf a}}{{\sf a}{\sf a}{\sf a}\overline{\sf a}}
-
\frac{{\sf c}\,d\overline{\sf a}}{
{\sf a}{\sf a}\overline{\sf a}\overline{\sf a}}
+
\frac{d{\sf c}}{{\sf a}{\sf a}\overline{\sf a}}.
\]

\medskip

In order to compute {\em at least partly}
the exterior derivative $d\sigma$ and to determine
the explicit expressions of {\em at least some of} the 
torsion coefficients $X_\smallbullet$ which will appear:
\[
\aligned
d\sigma
&
=
\big(2\,\alpha_1+\overline{\alpha}_1\big)
\wedge
\sigma
+
\\
&
\ \ \ \ \
+
X_1\,
\overline{\sigma}\wedge\sigma
+
X_2\,
\overline{\sigma}\wedge\rho
+
X_3\,
\overline{\sigma}\wedge\overline{\zeta}
+
X_4\,
\overline{\sigma}\wedge\zeta
\,+
\\
&
\ \ \ \ \ \ \ \ \ \ \ \ \ \ \ \ \ \ \ \ \ \ \ \ \ \ \
+
X_5\,
\sigma\wedge\rho
+
X_6\,
\sigma\wedge\overline{\zeta}
+
X_7\,
\sigma\wedge\zeta
\,+
\\
&
\ \ \ \ \ \ \ \ \ \ \ \ \ \ \ \ \ \ \ \ \ \ \ \ \ \ \
\ \ \ \ \ \ \ \ \ \ \ \ \ \ \ \ \ \ \ \ \ \ \ \ \ \ \
\ \ \ \ \ \ \ \ \ \ \ 
+
\rho\wedge\zeta,
\endaligned
\]
differentiate $\sigma = {\sf a}^2 \overline{\sf a}\, \sigma_0$, 
which gives:
\[
d\sigma
=
\big(
2\,{\sf a}\overline{\sf a}\,d{\sf a}
+
{\sf a}^2\,d\overline{\sf a}
\big)
\wedge
\sigma_0
+
{\sf a}^2\overline{\sf a}\,
d\sigma_0,
\]
obtain as an intermediate result:
\[
\aligned
d\sigma
&
=
\big(2\,\alpha^1+d\overline{\alpha}^1\big)
\wedge
\sigma
\,+
\\
&
\ \ \ \ \
+
{\sf a}^2\overline{\sf a}\,K\,\overline{\sigma}_0\wedge\sigma_0
+
{\sf a}^2\overline{\sf a}\,\overline{G}\,\overline{\sigma}_0
\wedge
\rho_0
+
{\sf a}^2\overline{\sf a}\,\overline{R}\,
\overline{\sigma}_0
\wedge
\overline{\zeta}_0
+
{\sf a}^2\overline{\sf a}\,B\,
\overline{\sigma}_0
\wedge
\zeta_0
\,+
\\
&
\ \ \ \ \ \ \ \ \ \ \ \ \ \ \ \ \ \ \ \ \ \ \ \ \ \ \ \ \ \ \ \ \ \ \ \
\ \ \ \ \ \ \ \ \ \ \ \ \ \ \ \ \ \ \ \ \ \ \ \ \ \ \ \ 
{\sf a}^2\overline{\sf a}\,F\,
\sigma_0
\wedge
\overline{\zeta}_0
+
{\sf a}^2\overline{\sf a}\,Q\,
\sigma_0
\wedge
\zeta_0
\,+
\\
&
\ \ \ \ \ \ \ \ \ \ \ \ \ \ \ \ \ \ \ \ \ \ \ \ \ \ \ \ \ \ \ \ \ \ \ \
\ \ \ \ \ \ \ \ \ \ \ \ \ \ \ \ \ \ \ \ \ \ \ \ \ \ \ \ \ \ \ \ \ \ \ \ 
\ \ \ \ \ \ \ \ \ \ \ \ \ \
+
{\sf a}^2\overline{\sf a}\,\,
\rho_0
\wedge
\zeta_0,
\endaligned
\]
and replace:
\[
\aligned
\sigma_0
&
=
\frac{1}{{\sf a}^2\overline{\sf a}}\,\sigma,
\\
\rho_0
&
=
-\,
\frac{\overline{\sf c}}{{\sf a}^2\overline{\sf a}^3}
-
\frac{{\sf c}}{{\sf a}^3\overline{\sf a}^2}\,\sigma
+
\frac{1}{{\sf a}\overline{\sf a}}\,\rho,
\\
\zeta_0
&
=
\frac{{\sf b}\overline{\sf c}-{\sf a}\overline{\sf a}\overline{\sf d}}{
{\sf a}^3\overline{\sf a}^3}
+
\frac{{\sf b}{\sf c}-{\sf a}\overline{\sf a}{\sf e}}{
{\sf a}^4\overline{\sf a}^2}\,
\sigma
-
\frac{{\sf b}}{{\sf a}^2\overline{\sf a}}\,\rho
+
\frac{1}{{\sf a}}\,\zeta,
\endaligned
\]
which, notably, provides:
\[
\aligned
X_2
&
=
\frac{1}{\overline{\sf a}^2}\,\overline{G}
-
\frac{{\sf b}}{{\sf a}\overline{\sf a}^2}\,B
-
\frac{\overline{\sf b}}{\overline{\sf a}^3}\,\overline{R}
+
\frac{\overline{\sf d}}{{\sf a}\overline{\sf a}^2},
\\
X_3
&
=
\frac{{\sf a}}{\overline{\sf a}^2}\,
\overline{R},
\\
X_4
&
=
\frac{1}{\overline{\sf a}}\,B
-
\frac{\overline{\sf c}}{{\sf a}\overline{\sf a}^2},
\\
X_6
&
=
\frac{1}{\overline{\sf a}}\,B,
\\
X_7
&
=
\frac{1}{{\sf a}}\,Q
-
\frac{{\sf c}}{{\sf a}^2\overline{\sf a}},
\endaligned
\]
other torsion coefficients being useless in what follows.

Proceeding similarly, the remaining two expressions of $d \rho$ and
of $d\zeta$: 
\[
\aligned
d\rho
&
=
\alpha^2\wedge\sigma
+
\overline{\alpha}^2\wedge\overline{\sigma}
+
\alpha^1\wedge\rho
+
\overline{\alpha}^1\wedge\rho
\,+
\\
&
\ \ \ \ \
+
Y_1\,\overline{\sigma}\wedge\sigma
+
Y_2\,\overline{\sigma}\wedge\rho
+
Y_3\,\overline{\sigma}\wedge\overline{\zeta}
+
Y_4\,\overline{\sigma}\wedge\zeta
\,+
\\
&
\ \ \ \ \ \ \ \ \ \ \ \ \ \ \ \ \ \ \ \ \ \ \ \,
+
\overline{Y}_2\,\sigma\wedge\rho
+
\overline{Y}_4\,\sigma\wedge\overline{\zeta}
+
\overline{Y}_3\,\sigma\wedge\zeta
\,+
\\
&
\ \ \ \ \ \ \ \ \ \ \ \ \ \ \ \ \ \ \ \ \ \ \ \ \ \ \ \ \ \ \ \ \ \ \ \
\ \ \ \ \ \ \ 
+
Y_8\,\rho\wedge\overline{\zeta}
+
\overline{Y}_8\,\rho\wedge\zeta
\,+
\\
&
\ \ \ \ \ \ \ \ \ \ \ \ \ \ \ \ \ \ \ \ \ \ \ \ \ \ \ \ \ \ \ \ \ \ \ \
\ \ \ \ \ \ \ \ \ \ \ \ \ \ \ \ \ \ \ \ \ \ \ \ \ 
+
\isqrt\,
\zeta\wedge\overline{\zeta},
\endaligned
\]
\[
\aligned
d\zeta
&
=
\alpha^3\wedge\sigma
+
\alpha^4\wedge\overline{\sigma}
+
\alpha^5\wedge\rho
+
\alpha^1\wedge\zeta
\,+
\\
&
\ \ \ \ \
+
Z_1\,\overline{\sigma}\wedge\sigma
+
Z_2\,\overline{\sigma}\wedge\rho
+
Z_3\,\overline{\sigma}\wedge\overline{\zeta}
+
Z_4\,\overline{\sigma}\wedge\zeta
\,+
\\
&
\ \ \ \ \ \ \ \ \ \ \ \ \ \ \ \ \ \ \ \ \ \ \ \,
+
Z_5\,\sigma\wedge\rho
+
Z_6\,\sigma\wedge\overline{\zeta}
+
Z_7\,\sigma\wedge\zeta
\,+
\\
&
\ \ \ \ \ \ \ \ \ \ \ \ \ \ \ \ \ \ \ \ \ \ \ \ \ \ \ \ \ \ \ \ \ \ \ \
\ \ \ \ \ \ \ 
+
Z_8\,\rho\wedge\overline{\zeta}
+
Z_9\,\rho\wedge\zeta
\,+
\\
&
\ \ \ \ \ \ \ \ \ \ \ \ \ \ \ \ \ \ \ \ \ \ \ \ \ \ \ \ \ \ \ \ \ \ \ \
\ \ \ \ \ \ \ \ \ \ \ \ \ \ \ \ \ \ \ \ \ \ \ \ \ 
+
Z_{10}\,
\zeta\wedge\overline{\zeta},
\endaligned
\]
can in principle be complemented by {\em explicit} expressions
of all the appearing torsion coefficients $Y_\smallbullet$ and
$Z_\smallbullet$, but only the single: 
\[
\overline{Y}_8
=
\frac{{\sf c}}{{\sf a}^2\overline{\sf a}}
+
\isqrt\,
\frac{\overline{\sf b}}{{\sf a}\overline{\sf a}}
\]
will be useful at this stage, because of absorption facts.

Indeed, following Cartan's procedure of determining the linear
subspace of torsion coefficients that are {\sl absorbable into
modified Maurer-Cartan forms} (\cite{ Olver-1995}), 
introduce modifications:
\[
\aligned
\widetilde{\alpha}^1
&
:=
\alpha^1
-
A_1\cdot\overline{\sigma}
-
B_1\cdot\sigma
-
C_1\cdot\rho
-
D_1\cdot\overline{\zeta}
-
E_1\cdot\zeta,
\\
\widetilde{\alpha}^2
&
:=
\alpha^2
-
A_2\cdot\overline{\sigma}
-
B_2\cdot\sigma
-
C_2\cdot\rho
-
D_2\cdot\overline{\zeta}
-
E_2\cdot\zeta,
\endaligned
\]
while the remaining three $1$-forms $\alpha^3$, $\alpha^4$, 
$\alpha^5$ are kept untouched, where $A_\smallbullet$,
$B_\smallbullet$, $C_\smallbullet$, $D_\smallbullet$,
$E_\smallbullet$ are functions defined on $M$, 
{\em i.e.} assumed to be {\em independent
of the group variables} ${\sf a}$, ${\sf b}$, ${\sf c}$, ${\sf d}$,
${\sf e}$. Replacing then $\alpha^1$ and $\alpha^2$ gives,
after collecting appropriately:
\[
\footnotesize
\aligned
d\sigma
&
=
\big(2\,\widetilde{\alpha}^1+\overline{\widetilde{\alpha}}^1\big)
\wedge
\sigma
\,+
\\
&
\ \ \ \ \
+
\overline{\sigma}\wedge\sigma\,
\big[
X_1-2\,A_1-\overline{B}_1
\big]
+
\overline{\sigma}\wedge\rho\,
\big[X_2\big]
+
\overline{\sigma}\wedge\overline{\zeta}\,
\big[X_3\big]
+
\overline{\sigma}\wedge\zeta\,
\big[X_4\big]
\,+
\\
&
\ \ \ \ \
+
\sigma\wedge\rho\,
\big[X_5-2\,C_1-\overline{C}_1\big]
+
\sigma\wedge\overline{\zeta}\,
\big[X_6-2\,D_1-\overline{E}_1\big]
+
\sigma\wedge\zeta\,
\big[X_7-2\,E_1-\overline{D}_1\big]
\,+
\\
&
\ \ \ \ \ \ \ \ \ \ \ \ \ \ \ \ \ \ \ \ \ \ \ \ \ \ \ \ \ \ \ \ \ \ \ \
\ \ \ \ \ \ \ \ \ \ \ \ \ \ \ \ \ \ \ \ \ \ \ \ \ \ \ \ \ \ \ \ \ \ \ \
\ \ \ \ \ \ \ \ \ \ \ \ \ \ \ \ \ \ \ \ \ \ \ \,
+
\rho\wedge\zeta,
\endaligned
\]
\[
\footnotesize
\aligned
d\rho
&
=
\widetilde{\alpha}^2\wedge\sigma
+
\overline{\widetilde{\alpha}}^2\wedge\overline{\sigma}
+
\widetilde{\alpha}^1\wedge\rho
+
\overline{\widetilde{\alpha}}^1\wedge\rho
\,+
\\
&
\ \ \ \ \
+
\overline{\sigma}\wedge\sigma\,
\big[Y_1+A_2-\overline{A}_2\big]
+
\overline{\sigma}\wedge\rho\,
\big[Y_2-\overline{C}_2+A_1+\overline{B}_1\big]
+
\overline{\sigma}\wedge\overline{\zeta}\,
\big[Y_3-\overline{E}_2\big]
+
\overline{\sigma}\wedge\zeta\,
\big[Y_4-\overline{D}_2\big]
\,+
\\
&
\ \ \ \ \
+
\sigma\wedge\rho\,
\big[\overline{Y}_2-C_2+\overline{A}_1+B_1\big]
+
\sigma\wedge\overline{\zeta}\,
\big[\overline{Y}_4-D_2\big]
+
\sigma\wedge\zeta\,
\big[\overline{Y}_3-E_2\big]
\,+
\\
&
\ \ \ \ \
+
\rho\wedge\overline{\zeta}\,
\big[Y_8-D_1-\overline{E}_1\big]
+
\rho\wedge\zeta\,
\big[\overline{Y}_8-\overline{D}_1-E_1\big]
\,+
\\
&
\ \ \ \ \
+
\isqrt\,
\zeta\wedge\overline{\zeta}.
\endaligned
\]
At first, the fact that in $d\rho$ the coefficient
$X_2$ of $\overline{ \sigma} \wedge \rho$, the
coefficient $X_3$ of $\overline{ \sigma} \wedge \overline{ \zeta}$
and the coefficient $X_4$ of $\overline{ \sigma} \wedge \zeta$
are left invariant, namely do not incorporate any $A_\smallbullet$,
$B_\smallbullet$, $C_\smallbullet$, $D_\smallbullet$, $E_\smallbullet$,
immediately explains the first steps of the following:

\begin{Lemma}
The three torsion coefficients $X_2$, $X_3$, $X_4$ are {\em essential}, 
and the same also holds true for
$\overline{ X}_6 + X_7 - 3\, \overline{ Y}_8$.
\end{Lemma}

Indeed, for this last, not immediately seen, linear combination, compute:
\[
\aligned
\overline{\widetilde{X}}_6
+
\widetilde{X}_7
-
3\,\overline{\widetilde{Y}}_8
&
=
\overline{X}_6
-
2\,\overline{D}_1
-
E_1
\,+
\\
&
+
X_7
-
2\,E_1
-
\overline{D}_1
\,+
\\
&
-\,
3\,\overline{Y}_8
+
3\,\overline{D}_1
+
3\,E_1
=
\overline{X}_6
+
X_7
-
3\,\overline{Y}_8.
\qed
\endaligned
\]

Observe here that when $R \not \equiv 0$, the essential torsion
coefficient $X_3 = \frac{ {\sf a}}{ \overline{\sf a}^2}\, \overline{
R}$ can be assigned the value $X_3 := 1$, which, after
relocalization to an open set on which $R$ is nowhere vanishing,
conducts to the normalization of the variable ${\sf a}$. This
observation led the first and the third authors in~\cite{
Merker-Sabzevari-III-1} to set up a natural {\em bifurcation} in
the concerned
equivalence problem, according to whether $R \equiv 0$ or $R
\not\equiv 0$. Several other potentially normalizable
essential torsion coefficients also appeared in the
advanced computational explorations performed in~\cite{
Merker-Sabzevari-III-1} which also concerned the diagonal
group parameter ${\sf a}$.

However, on the (simpler) way to construct
just a canonical Cartan connection, it is advisable
to ignore all such possible bifurcations, and to only look
at normalizations of group variables which come from $X_2$, 
$X_4$, $\overline{ X}_6 + X_7 - 3\, \overline{ Y}_8$,
{\em disregarding therefore $X_3$}.

In fact, setting equal to zero these three essential torsion 
coefficients provides, after elementary resolution, the
following three normalizations:
\[
\aligned
{\sf b}
&
=
{\sf a}\,
\big(
{\textstyle{\frac{\isqrt}{3}}}\,
\overline{Q}
-
\isqrt\,B
\big),
\\
{\sf c}
&
=
{\sf a}\overline{\sf a}\,\overline{B},
\\
{\sf d}
&
=
\overline{\sf a}\,
\big(
-\,G
-
{\textstyle{\frac{\isqrt}{3}}}\,Q
+
\isqrt\,\overline{B}
+
{\textstyle{\frac{\isqrt}{3}}}\,
\overline{Q}\,R
-
\isqrt\,B\,R
\big),
\endaligned
\]
which can be abbreviated as:
\[
{\sf b}
=
{\sf a}\,{\bf B}_0,
\ \ \ \ \ \ \ \ \ \ \ \ \ \ \ \ \ \ \ \ \ 
{\sf c}
=
{\sf a}\overline{\sf a}\,{\bf C}_0,\,
\ \ \ \ \ \ \ \ \ \ \ \ \ \ \ \ \ \ \ \ \ 
{\sf d}
=
\overline{\sf a}\,
{\bf D}_0,
\]
the three functions ${\bf B}_0$, ${\bf C}_0$, ${\bf D}_0$ being
functions defined on the basis $M$, independent of the group variables.

Then following Cartan and Olver
(\cite{ Olver-1995}), replace these group normalizations in:
\[
\aligned
\sigma
&
=
{\sf a}^2\overline{\sf a}\,\sigma_0,
\\
\rho
&
=
\overline{\sf c}\,\overline{\sigma}_0
+
{\sf c}\,\sigma_0
+
{\sf a}\overline{\sf a}\,\rho_0,
\\
\zeta
&
=
\overline{\sf d}\,\overline{\sigma}_0
+
{\sf e}\,\sigma_0
+
{\sf b}\,\rho_0
+
{\sf a}\,\zeta_0,
\endaligned
\]
which gives:
\[
\aligned
\sigma
&
=
{\sf a}^2\overline{\sf a}\,
\sigma_0
=:
{\sf a}^2\overline{\sf a}\,
\sigma_1,
\\
\rho
&
=
{\sf a}\overline{\sf a}\,
\big(
\overline{\bf C}_0\,\overline{\sigma}_0
+
{\bf C}_0\,\sigma_0
+
\rho_0
\big)
=:
{\sf a}\overline{\sf a}\,\rho_1,
\\
\zeta
&
=
{\sf e}\,\sigma_0
+
{\sf a}\,
\big(
\overline{\bf D}_0\,\overline{\sigma}_0
+
{\bf B}_0\,\rho_0
+
\zeta_0
\big)
=:
{\sf e}\,\sigma_1
+
{\sf a}\,\zeta_1,
\endaligned
\]
and restart the procedure of determining whether some group variables
are normalizable. In the present case, the second loop of Cartan's
procedure will happen to be the last one.

Of course, the new lifted coframe is:
\[
\left(\!\!
\begin{array}{c}
\overline{\sigma}
\\
\sigma
\\
\rho
\\
\overline{\zeta}
\\
\zeta
\end{array}
\!\!\right)
\,=\,
\left(\!
\begin{array}{ccccc}
{\sf a}\overline{\sf a}\overline{\sf a} & 0 & 0 & 0 & 0
\\
0 & {\sf a}{\sf a}\overline{\sf a} & 0 & 0 & 0
\\
0 & 0 & {\sf a}\overline{\sf a} & 0 & 0
\\
\overline{\sf e} & 0 & 0 & \overline{\sf a} & 0
\\
0 & {\sf e} & 0 & 0 & {\sf a}
\end{array}
\!\right)
\left(\!\!
\begin{array}{c}
\overline{\sigma}_1
\\
\sigma_1
\\
\rho_1
\\
\overline{\zeta}_1
\\
\zeta_1
\end{array}
\!\!\right),
\]
with two Maurer-Cartan $1$-forms (and their conjugates):
\[
\beta^1
:=
\frac{d{\sf a}}{{\sf a}},
\ \ \ \ \ \ \ \ \ \ \ \ \ \ \ \ \ \
\beta^2
:=
\frac{d{\sf e}}{{\sf a}^2\overline{\sf a}}
-
\frac{{\sf e}\,d{\sf a}}{{\sf a}^3\overline{\sf a}}.
\]
Some computations achieved in~\cite{ Merker-Sabzevari-III-1}
gave the explicit expressions of the new torsion coefficients in:
\[
\aligned
d\sigma
&
=
\big(
2\,\beta^1+\overline{\beta}^1
\big)
\wedge\sigma
\,+
\\
&
\ \ \ \ \
+
X_1'\,\overline{\sigma}\wedge\sigma
+
0
+
X_3'\,\overline{\sigma}\wedge\overline{\zeta}
+
0
\,+
\\
&
\ \ \ \ \
+
X_5'\,\sigma\wedge\rho
+
X_6'\,\sigma\wedge\overline{\zeta}
+
X_7'\,\sigma\wedge\zeta
+
\rho\wedge\zeta,
\endaligned
\]
the new $X_2' = 0$ and $X_4 ' = 0$ being zero thanks to the preceding
normalizations, and also the explicit expressions of the new torsion
coefficients appearing in:
\[
\aligned
d\rho
&
=
\big(\beta^1+\overline{\beta}^1\big)
\wedge\rho
\,+
\\
&
\ \ \ \ \
+
Y_1'\,\overline{\sigma}\wedge\sigma
+
Y_2'\,\overline{\sigma}\wedge\rho
+
Y_3'\,\overline{\sigma}\wedge\overline{\zeta}
+
Y_4'\,\overline{\sigma}\wedge\zeta
\,+
\\
&
\ \ \ \ \ 
+
\overline{Y}_2'\,\sigma\wedge\rho
+
\overline{Y}_4'\,\sigma\wedge\overline{\zeta}
+
\overline{Y}_3'\,\sigma\wedge\zeta
\,+
\\
&
\ \ \ \ \ 
+
\big(
{\textstyle{\frac{1}{3}}}\,
X_6'
+
{\textstyle{\frac{1}{3}}}\,
\overline{X}_7'
\big)\,
\rho\wedge\overline{\zeta}
+
\big(
{\textstyle{\frac{1}{3}}}\,
\overline{X}_6'
+
{\textstyle{\frac{1}{3}}}\,
X_7'
\big)\,
\rho\wedge\zeta
+
\isqrt\,
\zeta\wedge\overline{\zeta},
\endaligned
\]
but if just a Cartan connection is searched for as is admitted in the
present article, less computational efforts are demanded.

\begin{Proposition}
The appearing new torsion coefficient $Y_4'$ in $d\rho$ 
is essential, and assigned to $0$, it leads to a normalization of the last
non-diagonal group parameter ${\sf e}$ under the form:
\[
{\sf e}
=
{\sf a}\,
{\bf E}_0.
\]
\end{Proposition}

\proof
Indeed, the coefficient $Y_4'$ of the $2$-form $\overline{\sigma}
\wedge \zeta$ can visibly {\em not} be absorbed in $( \beta^1 + \overline{
\beta}^1 ) \wedge \rho$ by modifying:
\[
\widetilde{\beta}^1
:=
\beta^1
-
A_1\,\overline{\sigma}
-
B_1\,\sigma
-
C_1\,\rho
-
D_1\,\overline{\zeta}
-
E_1\,\zeta,
\]
which shows that $Y_4'$ is essential.  It
therefore only remains to explain how to shortly get at $Y_4'$.

Computing $d\rho_1$ and re-expressing it in terms of
the new coframe $\big\{ \overline{ \sigma}_1, \sigma_1, 
\rho_1, \overline{ \zeta}_1, \zeta_1 \big\}$ conducts to 
certain functions $T_{\smallbullet \smallbullet}^\smallbullet$
defined on $M$ that are independent of the group parameters:
\[
\aligned
d\rho_1
&
=
T_{\overline{\sigma}_1\sigma_1}^{\rho_1}\,
\overline{\sigma}_1\wedge\sigma_1
+
T_{\overline{\sigma}_1\rho_1}^{\rho_1}\,
\overline{\sigma}_1\wedge\rho_1
+
T_{\overline{\sigma}_1\overline{\zeta}_1}^{\rho_1}\,
\overline{\sigma}_1\wedge\overline{\zeta}_1
+
T_{\overline{\sigma}_1\zeta_1}^{\rho_1}\,
\overline{\sigma}_1\wedge\zeta_1
\,+
\\
&
\ \ \ \ \ \ \ \ \ \ \ \ \ \ \ \ \ \ \ \ \ \ \ \ \ \ \
+
T_{\sigma_1\rho_1}^{\rho_1}\,
\sigma_1\wedge\rho_1
+
T_{\sigma_1\overline{\zeta}_1}^{\rho_1}\,
\sigma_1\wedge\overline{\zeta}_1
+
T_{\sigma_1\zeta_1}^{\rho_1}\,
\sigma_1\wedge\zeta_1
\,+
\\
&
\ \ \ \ \ \ \ \ \ \ \ \ \ \ \ \ \ \ \ \ \ \ \ \ \ \ \ \ \ \ \ \ \ 
\ \ \ \ \ \ \ \ \ \ \ \ \ \ \ \ \ \ \ \
+
T_{\rho_1\overline{\zeta}_1}^{\rho_1}\,
\rho_1\wedge\overline{\zeta}_1
+
T_{\rho_1\zeta_1}^{\rho_1}\,
\rho_1\wedge\zeta_1
\,-
\\
&
\ \ \ \ \ \ \ \ \ \ \ \ \ \ \ \ \ \ \ \ \ \ \ \ \ \ \ \ \ \ \ \ \ 
\ \ \ \ \ \ \ \ \ \ \ \ \ \ \ \ \ \ \ \ \ \ \ \ \ \ \ \ \ \ \ \ \
\ \ \ \ \ \ \ \ \ \ \ \ \
-\,\isqrt\,
\overline{\zeta}_1\wedge\zeta_1.
\endaligned
\]
Differentiating $\rho = {\sf a} \overline{\sf a}\, \rho_1$ leads 
firstly to:
\[
\aligned
d\rho
&
=
\big(
{\textstyle{\frac{d{\sf a}}{{\sf a}}}}
+
{\textstyle{\frac{d\overline{\sf a}}{\overline{\sf a}}}}
\big)
\wedge\rho
\,+
\\
&
+
{\sf a}\overline{\sf a}\,
T_{\overline{\sigma}_1\sigma_1}^{\rho_1}\,
\overline{\sigma}_1\wedge\sigma_1
+
{\sf a}\overline{\sf a}\,
T_{\overline{\sigma}_1\rho_1}^{\rho_1}\,
\overline{\sigma}_1\wedge\rho_1
+
{\sf a}\overline{\sf a}\,
T_{\overline{\sigma}_1\overline{\zeta}_1}^{\rho_1}\,
\overline{\sigma}_1\wedge\overline{\zeta}_1
+
{\sf a}\overline{\sf a}\,
T_{\overline{\sigma}_1\zeta_1}^{\rho_1}\,
\underline{
\overline{\sigma}_1\wedge\zeta_1}
\,+
\\
&
\ \ \ \ \ \ \ \ \ \ \ \ \ \ \ \ \ \ \ \ \ \ \ \ \ \ \ \ \ \ \ 
+
{\sf a}\overline{\sf a}\,
T_{\sigma_1\rho_1}^{\rho_1}\,
\sigma_1\wedge\rho_1
+
{\sf a}\overline{\sf a}\,
T_{\sigma_1\overline{\zeta}_1}^{\rho_1}\,
\sigma_1\wedge\overline{\zeta}_1
+
{\sf a}\overline{\sf a}\,
T_{\sigma_1\zeta_1}^{\rho_1}\,
\sigma_1\wedge\zeta_1
\,+
\\
&
\ \ \ \ \ \ \ \ \ \ \ \ \ \ \ \ \ \ \ \ \ \ \ \ \ \ \ \ \ \ \ \ \ 
\ \ \ \ \ \ \ \ \ \ \ \ \ \ \ \ \ \ \ \ \ \ \ \ \ \ \ \
+
{\sf a}\overline{\sf a}\,
T_{\rho_1\overline{\zeta}_1}^{\rho_1}\,
\rho_1\wedge\overline{\zeta}_1
+
{\sf a}\overline{\sf a}\,
T_{\rho_1\zeta_1}^{\rho_1}\,
\rho_1\wedge\zeta_1
\,-
\\
&
\ \ \ \ \ \ \ \ \ \ \ \ \ \ \ \ \ \ \ \ \ \ \ \ \ \ \ \ \ \ \ \ \ 
\ \ \ \ \ \ \ \ \ \ \ \ \ \ \ \ \ \ \ \ \ \ \ \ \ \ \ \ \ \ \ \ \
\ \ \ \ \ \ \ \ \ \ \ \ \ \ \ \ \ \ \ \ \ \ \ \
-\,\isqrt\,
{\sf a}\overline{\sf a}\,
\underline{\overline{\zeta}_1\wedge\zeta_1},
\endaligned
\]
then replacing:
\[
\overline{\sigma}_1
=
\frac{1}{{\sf a}^2\overline{\sf a}}\,\sigma,
\ \ \ \ \ \ \ \ \ \ \ \ \ \ \ \ \ \ \ \ \ \ \ \
\zeta_1
=
-\,\frac{{\sf e}}{{\sf a}^3\overline{\sf a}}\,\sigma
+
\frac{1}{{\sf a}}\,\zeta,
\ \ \ \ \ \ \ \ \ \ \ \ \ \ \ \ \ \ \ \ \ \ \ \
\overline{\zeta}_1
=
-\,\frac{\overline{\sf e}}{{\sf a}\overline{\sf a}^3}\,\sigma
+
\frac{1}{\overline{\sf a}}\,\overline{\zeta},
\]
only the two underlined $2$-forms above 
contribute to $\overline{ \sigma} \wedge
\zeta$ in the final expression of $d\rho$, which gives:
\[
Y_4'
=
\frac{1}{{\sf a}\overline{\sf a}}\,
T_{\overline{\sigma}_1\zeta_1}^{\rho_1}
+
\isqrt\,
\frac{1}{{\sf a}\overline{\sf a}^2}\,
\overline{\sf e}.
\]
Setting $\overline{ Y}_4' = 0$ normalizes ${\sf e}$ under the form claimed.
\endproof

Performing therefore the obtained normalizations of all the nondiagonal
group parameters:
\[
{\sf b}
=
{\sf a}\,{\bf B}_0,
\ \ \ \ \ \ \ \ \ \ \ \ \ \ \
{\sf c}
=
{\sf a}\overline{\sf a}\,{\bf C}_0,
\ \ \ \ \ \ \ \ \ \ \ \ \ \ \
{\sf d}
=
\overline{\sf a}\,{\bf D}_0,
\ \ \ \ \ \ \ \ \ \ \ \ \ \ \
{\sf e}
=
{\sf a}\,{\bf E}_0,
\]
the reduced $G$-structure now involves only ${\sf a} \in \C^*$
and the new lifted coframe becomes:
\[
\aligned
\sigma
&
=
{\sf a}^2\overline{\sf a}\,\sigma_0
=:
{\sf a}^2\overline{\sf a}\,\sigma_2,
\\
\rho
&
=
{\sf a}\overline{\sf a}\,
\big(
\overline{\bf C}_0\,\overline{\sigma}_0
+
{\bf C}_0\,\sigma_0
+
\rho_0
\big)
=:
{\sf a}\overline{\sf a}\,
\rho_2
\\
\zeta
&
=
{\sf a}\,
\big(
\overline{\bf D}_0\,\overline{\sigma}_0
+
{\bf E}_0\,\sigma_0
+
{\bf B}_0\,\rho_0
+
\zeta_0
\big)
=:
{\sf a}\,\zeta_2,
\endaligned
\]
the single Maurer-Cartan $1$-form being (with its conjugate):
\[
\gamma^1
:=
\frac{d{\sf a}}{{\sf a}}.
\]
Inversion yields:
\[
\sigma_2
=
{\textstyle{\frac{1}{{\sf a}^2\overline{\sf a}}}}\,
\sigma,
\ \ \ \ \ \ \ \ \ \ \ \ \ \ \ \ \ \ \ \
\rho_2
=
{\textstyle{\frac{1}{{\sf a}\overline{\sf a}}}}\,
\rho,
\ \ \ \ \ \ \ \ \ \ \ \ \ \ \ \ \ \ \ \
\zeta_2
=
{\textstyle{\frac{1}{{\sf a}}}}\,
\zeta,
\]
whence immediately:
\[
\overline{\sigma}_2\wedge\sigma_2
=
{\textstyle{\frac{1}{{\sf a}^3\overline{\sf a}^3}}}\,
\overline{\sigma}\wedge\sigma,
\ \ \ \ \ \ \ \ \ 
\overline{\sigma}_2\wedge\rho_2
=
{\textstyle{\frac{1}{{\sf a}^2\overline{\sf a}^3}}}\,
\overline{\sigma}\wedge\rho,
\ \ \ \ \ \ \ \ \ 
\dots\dots,
\ \ \ \ \ \ \ \ \ 
\overline{\zeta}_2\wedge\zeta_2
=
{\textstyle{\frac{1}{{\sf a}\overline{\sf a}}}}\,
\overline{\zeta}\wedge\zeta.
\]
Then performing a last absorption thanks to the computer programs
of the second author:
\[
\lambda
:=
\frac{d{\sf a}}{{\sf a}}
+
\text{\rm linear combination of}\,
\big(
\overline{\sigma},\,\sigma,\,
\rho,
\overline{\zeta},\,\zeta
\big),
\]
the structure equations receive the final form:
\begin{equation}
\label{final-coframe}
\aligned
d\sigma
&
=
\big(2\,\lambda+\overline{\lambda}\big)
\wedge\sigma
+
\frac{{\sf a}}{\overline{\sf a}^2}\,\overline{R}\,\,
\overline{\sigma}\wedge\overline{\zeta}
+
\rho\wedge\zeta,
\\
d\rho
&
=
\big(\lambda+\overline{\lambda}\big)
\wedge\rho
+
\frac{1}{{\sf a}^2\overline{\sf a}^2}\,V_1\,
\overline{\sigma}\wedge\sigma
+
\frac{1}{\overline{\sf a}^2}\,V_3\,
\overline{\sigma}\wedge\overline{\zeta}
+
\frac{1}{{\sf a}^2\overline{\sf a}}\,\overline{V}_3\,
\sigma\wedge\rho
+
\isqrt\,\zeta\wedge\overline{\zeta},
\\
d\zeta
&
=
\lambda\wedge\zeta
+
\frac{1}{{\sf a}^2\overline{\sf a}^3}\,W_1\,
\overline{\sigma}\wedge\sigma
+
\frac{1}{{\sf a}\overline{\sf a}^3}\,W_2\,
\overline{\sigma}\wedge\rho
+
\frac{1}{\overline{\sf a}^3}\,W_3\,
\overline{\sigma}\wedge\overline{\zeta}
+
\frac{1}{{\sf a}\overline{\sf a}^2}\,W_4\,
\overline{\sigma}\wedge\zeta
\,+
\\
&
\ \ \ \ \ \ \ \ \ \ \ \ \ \ \ \ \ \ \ \ \ \ \ \ \ \ \ \ \ \ \ \ \
\ \ \ \ \ \ \ \ \ 
+
\frac{1}{{\sf a}^2\overline{\sf a}^2}\,W_5\,
\sigma\wedge\rho
+
\frac{1}{{\sf a}\overline{\sf a}^2}\,W_6\,
\sigma\wedge\overline{\zeta}
+
\frac{1}{{\sf a}^2\overline{\sf a}}\,W_7\,
\sigma\wedge\zeta
\,+
\\
&
\ \ \ \ \ \ \ \ \ \ \ \ \ \ \ \ \ \ \ \ \ \ \ \ \ \ \ \ \ \ \ \ \
\ \ \ \ \ \ \ \ \ \ \ \ \ \ \ \ \ \ \ \ \ \ \ \ \ \ \ \ \ \ \ \ \ 
\ \ \ \ \ 
+
\frac{1}{\overline{\sf a}^2}\,W_8\,\rho\wedge\overline{\zeta}
+
\frac{1}{{\sf a}\overline{\sf a}}\,W_9\,
\rho\wedge\zeta
\,+
\\
&
\ \ \ \ \ \ \ \ \ \ \ \ \ \ \ \ \ \ \ \ \ \ \ \ \ \ \ \ \ \ \ \ \
\ \ \ \ \ \ \ \ \ \ \ \ \ \ \ \ \ \ \ \ \ \ \ \ \ \ \ \ \ \ \ \ \ 
\ \ \ \ \ \ \ \ \ \ \ \ \ \ \ \ \ \ \ \ \ \ \ \ \ \ \ \ \ \
+
\frac{1}{\overline{\sf a}}\,W_{10}\,
\overline{\zeta}\wedge\zeta,
\endaligned
\end{equation}
in which all functions $\overline{ R}$, $V_\smallbullet$, 
$W_\smallbullet$ depend only on the variables of $M$, not
on the group parameter ${\sf a}$. 
But since $\frac{ d{\sf a}}{ {\sf a}}$ occuring in $\lambda$ is {\em closed},
$d\lambda$ is a linear combination of $2$-forms on $M$, namely
there are certain functions $I_{ \nu \mu}$ for
$\nu, \mu = \overline{ \sigma}, \sigma, \rho, \overline{ \zeta}, 
\zeta$ so that:
\begin{equation}
\label{lambda-end}
d\lambda
=
\sum_{\nu,\mu}\,
I_{\nu\mu}\,
\nu\wedge\mu,
\end{equation}
and Cartan's method then naturally stops (\cite{ Olver-1995}).

\begin{Theorem}
The $7$ differential $1$-forms $\overline{ \lambda}, 
\lambda, \overline{ \sigma}, \sigma, \rho, 
\overline{ \zeta}, \zeta$ define an absolute parallelism
on the principal bundle $P^7 := M \times \C^*$ satisfying the
structure equations~\thetag{ \ref{final-coframe}} 
and~\thetag{ \ref{lambda-end}} which reduces the 
local biholomorphic equivalence problem 
for Class $\text{\sf III}_{\text{\sf 1}}$ CR-generic submanifolds
$M^5 \subset \C^4$
with initial
structure group $G_{ \text{\sf III}_{\text{\sf 1}}}$ to an
$\{ e\}$-structure.\qed
\end{Theorem}

Lastly, elementary linear algebra shows
(\cite{ Merker-Sabzevari-III-1}, end of Section 12) 
that the Maurer-Cartan structure equations corresponding 
to the Lie algebra spanned by the above $7$ infinitesimal
CR automorphisms $S_2$, $S_1$, $T$, $L_2$, $L_1$, $D$, $R$
can be represented as:
\begin{equation}
\label{MC-model}
\aligned
d\alpha_{\sf c}
&
=
0,
\\
d\sigma_{\sf c}
&
=
\big(2\,\alpha_{\sf c}+\overline{\alpha}_{\sf c}\big)
\wedge
\sigma_{\sf c}
+
\rho_{\sf c}
\wedge
\zeta_{\sf c},
\\
d\rho_{\sf c}
&
=
\big(\alpha_{\sf c}+\overline{\alpha}_{\sf c}\big)
\wedge
\rho_{\sf c}
+
\isqrt\,
\zeta_{\sf c}
\wedge
\overline{\zeta}_{\sf c},
\\
d\zeta_{\sf c}
&
=
\alpha_{\sf c}
\wedge
\zeta_{\sf c},
\endaligned
\end{equation}
where $\alpha_{\sf c}$, $\sigma_{\sf c}$, $\zeta_{\sf c}$ are
$\C$-valued $1$-forms on the tangent bundle
$T M_{\sf c}$ to Beloshapka's cubic, and where $\rho_{\sf c}$ is
$\R$-valued $1$-form on $T M_{\sf c}$. 

Notably, when all invariants $\overline{ R}$, $V_\smallbullet$,
$W_\smallbullet$, $I_{\smallbullet \smallbullet}$, vanish in~\thetag{
\ref{final-coframe}} and in~\thetag{ \ref{lambda-end}}, renaming
$\lambda \mapsto \alpha_{\sf c}$ makes recover~\thetag{
\ref{MC-model}}.

\section{\bf Construction of a Cartan connection
(Proof of Theorem~\ref{main-theorem})}
\label{Cartan-connection}
\HEAD{\ref{Cartan-connection}.~Construction of a Cartan connection
(Proof of Theorem~\ref{main-theorem})}{
Jo\"el {\sc Merker}, Samuel {\sc Pocchiola}, Masoud {\sc Sabzevari}}

Introduce then the set of vector fields:
\[
\big\{
{\sf e}_{\overline{\alpha}},\,
{\sf e}_\alpha,\,
{\sf e}_{\overline{\sigma}},\,
{\sf e}_\sigma,\,
{\sf e}_\rho,\,
{\sf e}_{\overline{\zeta}},\,
{\sf e}_\zeta
\big\}
\]
that is dual to 
$\big\{ \overline{ \alpha}_{\sf c}, \, 
\alpha_{\sf c},\, \overline{ \sigma}_{\sf c},\,
\sigma_{\sf c},\, \rho_{\sf c},\, 
\overline{ \zeta}_{\sf c},\, \zeta_{\sf c} \big\}$, 
hence defines the structure of a $7$-dimensional real Lie algebra:
\[
\mathfrak{g}^7
\cong
\mathfrak{aut}_{\text{\sc cr}}(M_{\sf c}^5).
\]
From~\thetag{ \ref{MC-model}}, the Lie bracket structure of
$\mathfrak{ g}^7$ is:
\[
\aligned
\big[
{\sf e}_{\overline{\alpha}},\,
{\sf e}_{\overline{\sigma}}
\big]
&
=
-\,2\,{\sf e}_{\overline{\sigma}},
\\
\big[
{\sf e}_{\overline{\alpha}},\,
{\sf e}_\sigma
\big]
&
=
-\,{\sf e}_\sigma,
\\
\big[
{\sf e}_{\overline{\alpha}},\,
{\sf e}_\rho
\big]
&
=
-\,{\sf e}_\rho,
\\
\big[
{\sf e}_{\overline{\alpha}},\,
{\sf e}_{\overline{\zeta}}
\big]
&
=
-\,
{\sf e}_{\overline{\zeta}},
\endaligned
\ \ \ \ \ \ \ \ \ \ \
\aligned
\big[
{\sf e}_\alpha,\,
{\sf e}_{\overline{\sigma}}
\big]
&
=
-\,{\sf e}_{\overline{\sigma}},
\\
\big[
{\sf e}_\alpha,\,
{\sf e}_\sigma
\big]
&
=
-\,2\,
{\sf e}_\sigma,
\\
\big[
{\sf e}_\alpha,\,
{\sf e}_\rho
\big]
&
=
-\,{\sf e}_\rho,
\\
\big[
{\sf e}_\alpha,\,
{\sf e}_\zeta
\big]
&
=
-\,{\sf e}_\zeta,
\endaligned
\ \ \ \ \ \ \ \ \ \ \
\aligned
\big[
{\sf e}_\rho,\,
{\sf e}_{\overline{\zeta}}
\big]
&
=
-\,{\sf e}_{\overline{\sigma}},
\\
\big[
{\sf e}_\rho,\,
{\sf e}_\sigma
\big]
&
=
-\,{\sf e}_\sigma,
\endaligned
\ \ \ \ \ \ \ \ \ \ \
\aligned
\big[
{\sf e}_{\overline{\zeta}},\,
{\sf e}_\zeta
\big]
&
=
-\,\isqrt\,
{\sf e}_\rho,
\endaligned
\]
unwritten brackets being zero.

Next, let $\mathfrak{ h} \cong \R^2$ be the real Lie algebra spanned
by $\{ {\sf e}_{\overline{ \alpha}},\, {\sf e}_\alpha \}$.
Let $P^7$ be $M \times \C^*$ equipped with some local coordinates 
on $M$ and with the fiber coordinates $( {\sf a}, \overline{\sf a})$.

\begin{Theorem}
In terms of the $7$ differential $1$-forms $\overline{ \lambda}$,
$\lambda$, $\overline{ \sigma}$, $\sigma$, $\rho$, $\overline{
\zeta}$, $\zeta$, obtained by reducing to an $\{ e\}$-structure 
the biholomorphic equivalence problem for Class $\text{\sf
III}_{\text{\sf 1}}$ CR-generic submanifolds $M^5 \subset \C^4$,
the $1$-form $\varpi$ with value in $\mathfrak{ g}^7$ defined at 
an arbitrary point $p \in P$ for every tangent vector 
$v_p \in T_p P^7$ by:
\[
\varpi_p(v_p)
:=
\overline{\lambda}_p(v_p)
\cdot
{\sf e}_{\overline{\alpha}}
+
\lambda_p(v_p)
\cdot
{\sf e}_{\alpha}
+
\overline{\sigma}_p(v_p)
\cdot
{\sf e}_{\overline{\sigma}}
+
\sigma_p(v_p)
\cdot
{\sf e}_\sigma
+
\rho_p(v_p)
\cdot
{\sf e}_\rho
+
\overline{\zeta}_p(v_p)
\cdot
{\sf e}_{\overline{\zeta}}
+
\zeta_p(v_p)
\cdot
{\sf e}_\zeta
\]
defines a canonical Cartan connection on $P^7 = M^5 \times \C^*$.
\end{Theorem}

\proof
Since by construction $\overline{ \lambda}, \lambda, \overline{
\sigma}, \sigma, \rho, \overline{\zeta}, \zeta$ span the cotangent
bundle to $P^7$ at every point, condition {\bf (i)} for a Cartan
connection that $\varpi_p \colon T_p P^7 \longrightarrow \mathfrak{
g}^7$ be an isomorphism at every point
$p \in P$ is automatically satisfied.

As well, the condition {\bf (ii)} that $\varpi (Y^\dag) = {\sf y}$
comes directly from the fact that: 
\[
\lambda
=
\frac{d{\sf a}}{{\sf a}}
+
\text{\rm linear combinations of}\,
\big(
\overline{\sigma},\sigma,\rho,\overline{\zeta},\zeta
\big).
\]

Lastly, condition {\bf (iii)} is known 
(\cite{ Sternberg-1964})
to be equivalent to its infinitesimal
counterpart:
\[
{\sf Lie}_{{\sf e}_\alpha^\dag}
\big(\varpi)
\,=\,
-\,{\sf ad}_{{\sf e}_\alpha}
\circ
\varpi, 
\ \ \ \ \ \ \ \ \ \ \ \ \ \ \ \ \ \
\text{\rm and}
\ \ \ \ \ \ \ \ \ \ \ \ \ \ \ \ \ \
{\sf Lie}_{{\sf e}_{\overline{\alpha}}^\dag}
\big(\varpi)
\,=\,
-\,{\sf ad}_{{\sf e}_{\overline{\alpha}}}
\circ
\varpi, 
\]
where ${\sf Lie}_X (\varpi)$ denotes the Lie derivative of $\varpi$
with respect to a vector field $X$ on $P$, and where ${\sf ad}_{\sf
k}$ is the linear map $\mathfrak{ k} \longrightarrow \mathfrak{ k}$
defined on a Lie algebra $\mathfrak{ k}$ by ${\sf ad}_{\sf k} ( {\sf
l}) := [ {\sf k}, {\sf l} ]$.

Now, condition {\bf (ii)} verified at the moment shows\,\,---\,\,
using $\lrcorner$ to denote interior product\,\,---\,\,that in Cartan's
formula:
\[
{\sf Lie}_{{\sf e}_\alpha^\dag}
(\varpi)
=
{\sf e}_{\alpha^\dag}
\lrcorner\,
d\varpi
+
\zero{d\big(
{\sf e}_\alpha^\dag
\lrcorner\,\varpi
\big)}
\ \ \ \ \ \ \ \ \ \ \ \ \ \ \ \ \ \
\text{\rm and in:}
\ \ \ \ \ \ \ \ \ \ \ \ \ \ \ \ \ \
{\sf Lie}_{{\sf e}_{\overline{\alpha}}^\dag}
(\varpi)
=
{\sf e}_{\overline{\alpha}^\dag}
\lrcorner\,
d\varpi
+
\zero{d\big(
{\sf e}_{\overline{\alpha}}^\dag
\lrcorner\,\varpi
\big)},
\]
the second terms, differentiating a constant, drop, so that
verifying {\bf (iii)} amounts to checking the two coincidences:
\[
{\sf e}_{\alpha^\dag}
\lrcorner\,
d\varpi
\,=\,
-\,{\sf ad}_{{\sf e}_\alpha}
\big(\varpi), 
\ \ \ \ \ \ \ \ \ \ \ \ \ \ \ \ \ \
\text{\rm and}
\ \ \ \ \ \ \ \ \ \ \ \ \ \ \ \ \ \
{\sf e}_{\overline{\alpha}^\dag}
\lrcorner\,
d\varpi
\,=\,
-\,{\sf ad}_{{\sf e}_{\overline{\alpha}}}
\big(\varpi). 
\]

But from the structure equations~\thetag{ \ref{final-coframe}} 
and~\thetag{ \ref{lambda-end}} satisfied by the $7$ two-forms 
occuring in $d\varpi$, it is clear that:
\[
\aligned
{\sf e}_\alpha^\dag
\lrcorner\,
d\overline{\lambda}
=
0,
\ \ \ \ \ \ \ \ \ \
{\sf e}_\alpha^\dag
\lrcorner\,
d\lambda
&
=
0,
\ \ \ \ \ \ \ \ \ \
{\sf e}_\alpha^\dag
\lrcorner\,
d\overline{\sigma}
=
\overline{\sigma},
\ \ \ \ \ \ \ \ \ \
{\sf e}_\alpha^\dag
\lrcorner\,
d\sigma
=
2\,\sigma,
\\
{\sf e}_\alpha^\dag
\lrcorner\,
d\rho
&
=
\rho,
\ \ \ \ \ \ \ \ \ \
{\sf e}_\alpha^\dag
\lrcorner\,
d\overline{\zeta}
=
0,
\ \ \ \ \ \ \ \ \ \ \
{\sf e}_\alpha^\dag
\lrcorner\,
d\zeta
=
\zeta,
\endaligned
\]
while the Lie bracket structure shows on the other hand that:
\[
\footnotesize
\aligned
-\,{\sf ad}_{{\sf e}_\alpha}
\circ
\varpi
&
=
-\,{\sf ad}_{{\sf e}_\alpha}
\circ
\Big(
\overline{\lambda}\cdot{\sf e}_{\overline{\alpha}}
+
\lambda\cdot{\sf e}_\alpha
+
\overline{\sigma}\cdot{\sf e}_{\overline{\sigma}}
+
\sigma\cdot{\sf e}_\sigma
+
\rho\cdot{\sf e}_\rho
+
\overline{\zeta}\cdot{\sf e}_{\overline{\zeta}}
+
\zeta\cdot{\sf e}_\zeta
\Big)
\\
&
=
-\,\overline{\lambda}\,
\big[{\sf e}_\alpha,\,{\sf e}_{\overline{\alpha}}\big]
-
\lambda\,\big[{\sf e}_\alpha,\,{\sf e}_\alpha\big]
-
\overline{\sigma}\,\big[{\sf e}_\alpha,\,{\sf e}_{\overline{\sigma}}\big]
-
\sigma\,\big[{\sf e}_\alpha,\,{\sf e}_\sigma\big]
-
\rho\,\big[{\sf e}_\alpha,\,{\sf e}_\rho\big]
-
\overline{\zeta}\,\big[{\sf e}_\alpha,\,{\sf e}_{\overline{\zeta}}\big]
-
\zeta\,\big[{\sf e}_\alpha,\,{\sf e}_\zeta\big]
\\
&
=
0
+
0
+
\overline{\sigma}\cdot{\sf e}_{\overline{\sigma}}
+
2\,\sigma\cdot{\sf e}_{\sigma}
+
\rho\cdot{\sf e}_\rho
+
0
+
\zeta\cdot{\sf e}_\zeta,
\endaligned
\]
so that the first coincidence holds; the second one is treated similarly.
\endproof

\section{\bf Beyond Cartan connections}
\label{beyond-connections}
\HEAD{\ref{beyond-connections}.~Beyond Cartan connections}{
Jo\"el {\sc Merker}, Samuel {\sc Pocchiola}, Masoud {\sc Sabzevari}}

This article announces that, beyond plain linear algebra
considerations that are sufficient to construct $\{ e\}$-structures
and Cartan connections associated to the local biholomorphic
equivalence problem for CR-generic $\mathcal{ C}^\omega$ submanifolds
$M^5 \subset \C^4$ belonging to Class $\text{\sf III}_{\text{\sf 1}}$,
the memoir~\cite{ Merker-Sabzevari-III-1} presents {\em explicit}
expressions of the incoming invariants $\overline{ R}$,
$V_\smallbullet$, $W_\smallbullet$, $I_{\smallbullet \smallbullet}$, 
it examines the {\em nonlinear} relations these invariants may share,
and it conducts a ramified analysis of the {\em bifurcation tree} of
possible further normalizations for the diagonal last group
parameter ${\sf a}$ which entails non-uniqueness of the concerned
Cartan connections or $\{ e\}$-structures.

A unified presentation of recent results on Cartan
equivalences for all six general classes $\text{\sf I}$, $\text{\sf
II}$, $\text{\sf III}_{\text{\sf 1}}$, $\text{\sf III}_{\text{\sf
2}}$, $\text{\sf IV}_{\text{\sf 1}}$, $\text{\sf IV}_{\text{\sf
2}}$ of CR manifolds of dimension
$\leqslant 5$ is upcoming (\cite{
Merker-Pocchiola-Sabzevari-5-CR-I}).

\vfill\end{document}